\documentclass{article}
\usepackage{amsmath,amssymb}
\newtheorem{df}{Definition}[section]
\newtheorem{thm}[df]{Theorem}

\newtheorem{cor}[df]{Corollary}
\title{Frames of directional wavelets\\ on $n$-dimensional spheres}
\author{Ilona Iglewska-Nowak\footnote{West Pomeranian University of Technology, School of Mathematics, al. Piast\'ow 17, 70--310 Szczecin, Poland}}

\begin{document}

\maketitle

\bibliographystyle{amsplain}

\begin{abstract}The major goal of the paper is to prove that discrete frames of (directional) wavelets derived from an approximate identity exist. Additionally, a kind of energy conservation property is shown to hold in the case when a wavelet family is not its own reconstruction family. Although an additional constraint on the spectrum of the wavelet family must be satisfied, it is shown that all the wavelets so far defined in the literature possess this property.
\end{abstract}

\begin{bfseries}Keywords:\end{bfseries} spherical wavelets, approximate identities, frames, $n$-spheres\\
\begin{bfseries}AMS Classification:\end{bfseries} 42C40

\section{Introduction}

The present paper is a continuation of the series \cite{IIN14CWT,IIN14PW,IIN15WF,IIN15DW}. A wide class of zonal wavelets is constructed, which contains all the so far studied wavelet families (derived from approximate identities). For the case of non--zonal wavelets, the definition of the wavelet transform given in~\cite{IIN14CWT} is weakened, which is a price for having fully discrete directional wavelet frames. The frame constructions on $n$--dimensional spheres known to the author are all based on zonal wavelets \cite{GM09b,IIN15WF}, therefore, the present paper seems to be the first approach to discretize a directional wavelet transform.

The paper is organized as follows. In Section~\ref{sec:sphere} we recapitulate basic facts about spherical functions, frames, and wavelets derived from approximate identities. We also present a wide class of functions that can serve as wavelet families. In Section~\ref{sec:semi_cont_frames} we show that the wavelet transform with respect to such wavelet families can be discretized with respect to the scale parameter, and in Section~\ref{sec:discrete_frames} we prove the existence of fully discrete wavelet frames.

\section{Preliminaries}\label{sec:sphere}

\subsection{Functions on the sphere}

By $\mathcal{S}^n$ we denote the $n$--dimensional unit sphere in $n+1$--dimensional Euclidean space~$\mathbb{R}^{n+1}$ with the rotation--invariant measure~$d\sigma_n$ normalized such that
$$
\Sigma_n=\int_{\mathcal{S}^n}d\sigma_n=\frac{2\pi^{(n+1)/2}}{\Gamma\bigl((n+1)/2\bigr)}.
$$
The surface element $d\sigma_n$ is explicitly given by
$$
d\sigma_n=\sin^{n-1}\theta_1\,\sin^{n-2}\theta_2\dots\sin\theta_{n-1}d\theta_1\,d\theta_2\dots d\theta_{n-1}d\varphi,
$$
where $(\theta_1,\theta_2,\dots,\theta_{n-1},\varphi)\in[0,\pi]^{n-1}\times[0,2\pi)$ are spherical coordinates satisfying
\begin{align*}
x_1&=\cos\theta_1,\\
x_2&=\sin\theta_1\cos\theta_2,\\
x_3&=\sin\theta_1\sin\theta_2\cos\theta_3,\\
&\dots\\
x_{n-1}&=\sin\theta_1\sin\theta_2\dots\sin\theta_{n-2}\cos\theta_{n-1},\\
x_n&=\sin\theta_1\sin\theta_2\dots\sin\theta_{n-2}\sin\theta_{n-1}\cos\varphi,\\
x_{n+1}&=\sin\theta_1\sin\theta_2\dots\sin\theta_{n-2}\sin\theta_{n-1}\sin\varphi.
\end{align*}
$\left<x,y\right>$ or $x\cdot y$ stand for the scalar product of vectors with origin in~$O$ and an endpoint on the sphere. As long as it does not lead to misunderstandings, we identify these vectors with points on the sphere. By their geodesic distance we mean $\angle(x,y)=\arccos(x,y)$.

A function is called zonal if its value depends only on $\theta=\theta_1=\left<\hat e,x\right>$, where~$\hat e$ is the north pole of the sphere
$$
\hat e=(1,0,0,\dots,0).
$$
It is invariant with respect to the rotation about the axis through~$O$ and~$\hat e$. We identify zonal functions with functions over the interval $[-1,1]$, i.e., whenever it does not lead to mistakes, we write
$$
f(x)=f(\cos\theta_1).
$$
 The subspace of $p$--integrable zonal functions is isomorphic to and will be identified with the space~$\mathcal L_\lambda^p$, and the norms satisfy
$$
\|f\|_{\mathcal{L}^p(\mathcal{S}^n)}=\|f\|_{\mathcal{L}_\lambda^p([-1,1])}
$$
for
$$
\|f\|_{\mathcal{L}^p(\mathcal{S}^n)}:=\left[\frac{1}{\Sigma_n}\int_{\mathcal S^n}|f(x)|^p\,d\sigma_n(x)\right]^{1/p}
$$
and
$$
\lambda:=\frac{n-1}{2}.
$$

The scalar product of $f,g\in\mathcal L^2(\mathcal S^n)$ is defined by
$$
\left<f,g\right>_{\mathcal L^2(\mathcal S^n)}=\frac{1}{\Sigma_n}\int_{\mathcal S^n}\overline{f(x)}\,g(x)\,d\sigma_n(x),
$$
such that $\|f\|_2^2=\left<f,f\right>$.

Gegenbauer polynomials $C_l^\lambda$ of order~$\lambda\in\mathbb R$ and degree $l\in\mathbb{N}_0$ are defined in terms of their generating function
$$
\sum_{l=0}^\infty C_l^\lambda(t)\,r^l=\frac{1}{(1-2tr+r^2)^\lambda},\qquad t\in[-1,1].
$$
A set of Gegenbauer polynomials $\bigl\{C_l^\lambda\bigr\}_{l\in\mathbb N_0}$ builds a complete orthogonal system on $[-1,1]$ with weight $(1-t^2)^{\lambda-1/2}$. Consequently, it is an orthogonal basis for zonal functions on the $2\lambda+1$--dimensional sphere.

Let $Q_l$ denote a polynomial on~$\mathbb{R}^{n+1}$ homogeneous of degree~$l$, i.e., such that $Q_l(az)=a^lQ_l(z)$ for all $a\in\mathbb R$ and $z\in\mathbb R^{n+1}$, and harmonic in~$\mathbb{R}^{n+1}$, i.e., satisfying $\nabla^2Q_l(z)=0$, then $Y_l(x)=Q_l(x)$, $x\in\mathcal S^n$, is called a hyperspherical harmonic of degree~$l$. The set of hyperspherical harmonics of degree~$l$ restricted to~$\mathcal S^n$ is denoted by $\mathcal H_l(\mathcal S^n)$. Hyperspherical harmonics of distinct degrees are orthogonal to each other. The number of linearly independent hyperspherical harmonics of degree~$l$ is equal to
$$
N=N(n,l)=\frac{(n+2l-1)(n+l-2)!}{(n-1)!\,l!}.
$$

In this paper, we will be working with the orthogonal basis for~$\mathcal L^2(\mathcal S^n)=\overline{\bigoplus_{l=0}^\infty\mathcal H_l}$, consisting of hyperspherical harmonics given by
$$
Y_l^k(x)=A_l^k\prod_{\tau=1}^{n-1}C_{k_{\tau-1}-k_\tau}^{\frac{n-\tau}{2}+k_\tau}(\cos\theta_\tau)\sin^{k_\tau}\!\theta_\tau\cdot e^{\pm ik_{n-1}\varphi}
$$
with $l=k_0\geq k_1\geq\dots\geq k_{n-1}\geq0$, $k$ being a sequence $(k_1,\dots,\pm k_{n-1})$ of integer numbers, and normalization constants~$A_l^k$. The set of non-increasing sequences~$k$ in $\mathbb N_0^{n-1}\times\mathbb Z$ with elements bounded by~$l$ will be denoted by $\mathcal M_{n-1}(l)$.

Funk--Hecke theorem states that for $f\in\mathcal L_\lambda^1\bigl([-1,1])$ and $Y_l\in\mathcal H_l(\mathcal S^n)$, $l\in\mathbb{N}_0$,
\begin{equation}\label{eq:FH}\begin{split}
\int_{\mathcal S^n}&Y_l(y)\,f(x\cdot y)\,d\sigma_n(y)\\
&=Y_l(x)\cdot\frac{(4\pi)^\lambda\,l!\,\Gamma(\lambda)}{(2\lambda+l-1)!}\int_{-1}^1 f(t)\,C_l^\lambda(t)\left(1-t^2\right)^{\lambda-1/2}dt.
\end{split}\end{equation}

Every $\mathcal{L}^1(\mathcal S^n)$--function~$f$ can be expanded into Laplace series of hyperspherical harmonics by
$$
f\sim\sum_{l=0}^\infty f_l,
$$
where $f_l$ is given by
\begin{equation*}
f_l(x)=\frac{\lambda+l}{\lambda\Sigma_n}\int_{\mathcal S^n}f(y)\,C_l^\lambda(x\cdot y)\,d\sigma_n(y).
\end{equation*}
For zonal functions we obtain by~\eqref{eq:FH} the representation
$$
f_l(t)=\widehat f(l)\,C_l^\lambda(t),\qquad t=\cos\theta,
$$
with Gegenbauer coefficients
$$
\widehat f(l)=c(l,\lambda)\int_{-1}^1 f(t)\,C_l^\lambda(t)\left(1-t^2\right)^{\lambda-1/2}dt
$$
for some constants $c(l,\lambda)$. The series
$$
\sum_{l=0}^\infty\widehat f(l)\,C_l^\lambda(t)
$$
is called Gegenbauer expansion of~$f$.

For $f,h\in\mathcal L^1(\mathcal S^n)$, $h$ zonal, their convolution $f\ast h$ is defined by
\begin{equation*}
(f\ast h)(x)=\frac{1}{\Sigma_n}\int_{\mathcal S^n}f(y)\,h(x\cdot y)\,d\sigma_n(y).
\end{equation*}
With this notation we have
$$
f_l(x)=\frac{\lambda+l}{\lambda}\,\bigl(f\ast C_l^\lambda\bigr)(x),
$$
hence, the function $\mathcal K_l^\lambda:=\frac{\lambda+l}{\lambda}\,C_l^\lambda$ is the reproducing kernel for~$\mathcal H_l(\mathcal S^n)$, and Funk--Hecke formula can be written as
$$
Y_l\ast f=\frac{\lambda}{\lambda+l}\,\widehat f(l)\,Y_l.
$$

Further, any function $f\in\mathcal L^2(\mathcal S^n)$ has a unique representation as a mean--convergent series
$$
f(x)=\sum_{l=0}^\infty\sum_{k\in\mathcal M_{n-1}(l)} a_l^k\,Y_l^k(x),\qquad x\in\mathcal S^n,
$$
where
$$
a_l^k=a_l^k(f)=\frac{1}{\Sigma_n}\int_{\mathcal S^n}\overline{Y_l^k(x)}\,f(x)\,d\sigma_n(x)=\left<Y_l^k,f\right>.
$$
We call $a_l^k$ the Fourier coefficients of the function~$f$. Convolution with a zonal function can be then written as
$$
f\ast g=\sum_{l=0}^\infty\sum_{k\in\mathcal M_{n-1}(l)} \frac{\lambda}{\lambda+l}\,a_l^k(f)\,\widehat g(l)\,Y_l^k
$$
and for zonal functions the following relation
\begin{equation}\label{eq:hatfl_vs_al0f}
\widehat f(l)=A_l^0\cdot a_l^0(f)
\end{equation}
between Fourier and Gegenbauer coefficients holds.

The set of rotations of~$\mathbb R^{n+1}$ is denoted by $SO(n+1)$. It is isomorphic to the set of square matrices of degree $n+1$ with determinant~$1$. The $n$--dimensional sphere can be identified with the class of left cosets of $SO(n+1)$ mod $SO(n)$,
$$
\mathcal S^n=SO(n+1)/SO(n),
$$
cf.~\cite[Sec.~I.2]{Vilenkin}. Consequently, $SO(n+1)$ can be interpreted as the Cartesian product of~$\mathcal S^J$, $J=1,\,2,\dots,\,n$. The invariant unit measure on $SO(n+1)$ will be denoted by~$\nu$ and is given by
$$
d\nu=\left(\prod_{j=1}^n\Sigma_n\right)^{-1}\prod d\sigma_n
$$
for $\Upsilon\in SO(n+1)$ parameterized by its Euler angles
\begin{equation*}\begin{split}
&\theta_1^{(n)},\,\theta_2^{(n)},\,\dots,\,\theta_{n-1}^{(n)},\,\varphi^{(n)},\\
&\theta_1^{(n-1)},\,\theta_2^{(n-1)},\,\dots,\,\theta_{n-2}^{(n-1)},\,\varphi^{(n-1)},\\
&\dots,\\
&\theta_1^{(2)},\,\varphi^{(2)},\\
&\varphi^{(1)},
\end{split}\end{equation*}
compare \cite[Sec.~IX.1.3--4]{Vilenkin}.

For further details on this topic, compare the textbooks~\cite{Vilenkin} and~\cite{AH12}.

\subsection{Frames in reproducing kernel Hilbert spaces}

The goal of this paper is to prove that discrete spherical wavelet frames exist. For this, we need the following definition.

\begin{df} A family of vectors $\{g_x,\,x\in X\}\subset\mathcal{H}$ in a Hilbert space $\mathcal{H}$ indexed by a measure space $X$ with a positive measure $\mu$
is called a~\emph{frame with weight}~$\mu$ if the mapping $x\mapsto g_x$ is weakly measurable, i.e., $x\mapsto\left<g_x,u\right>$ is measurable, and if
for some $0\leq \epsilon<1$ we have
\begin{equation}  \label{eq:framedef}
(1-\epsilon)\,\|u\|^2\leq \int_{X}|\left<g_x,u\right>|^2\,d\mu(x) \leq(1+\epsilon)\,\|u\|^2.
\end{equation}
for all $u\in\mathcal{H}$. Equivalently, the frame condition reads
\begin{equation*}
\left| \int_X |\left<g_x,u\right>|^2 \,d\mu(x) - \|u\|^2 \right| \leq\epsilon\, \|u\|^2.
\end{equation*}
If $\epsilon=0$, we call it a \emph{tight frame}. The numbers $A:=1-\epsilon$ and $B:=1+\epsilon$ are called \emph{frame bounds}.
\end{df}

\begin{bfseries}Remark.\end{bfseries} Alternatively, one can require
\begin{equation}\label{eq:framedef_AB}
A\,\|u\|^2\leq \int_{X}|\left<g_x,u\right>|^2\,d\mu(x) \leq B\,\|u\|^2.
\end{equation}
for some $0<A\leq B$. A multiplication of this inequality by~$\frac{2}{A+B}$ leads to~\eqref{eq:framedef} with the weight function modified by this factor. The choice of~\eqref{eq:framedef} as a definition simplifies the formulation of the next theorem.

\begin{thm}\label{thm:frames_discretization}
Let $\{g_x,\,x\in X\}\subset\mathcal{H}$ be a frame in~$\mathcal H$ with weight~$\mu$ and~$\epsilon$ as in~\eqref{eq:framedef}. Further, suppose $\Lambda\subset X$ and $\lambda$ is a measure on~$\Lambda$. Then the family of functions $\{g_y:\,y\in\Lambda\}\subset\mathcal L^2(X,d\mu)$ is a frame with weight~$\lambda$ for the image of the mapping
$$
\mathcal{F} : \mathcal{H} \to \mathcal L^2(X,d\mu),\qquad\mathcal{F}u(x)=\left<g_x,u\right>
$$
if
$$
\left| \int_\Lambda |\left<g_y,u\right>|^2 \,d\lambda(y)- \int_X |\left<g_x,u\right>|^2 \,d\mu(x)\right| \leq\delta\, \|u\|^2
$$
for some $\delta<1-\epsilon$.
\end{thm}

More details on this topic can be found in~\cite{IH10}, \cite{CJ99}, and~\cite{oC03}.

\subsection{Continuous wavelet transform -- bilinear theory}\label{subs:bilinearwv}

Nonzonal bilinear wavelets derived from an approximate identity were introduced in~\cite{EBCK09} and further investigated in~\cite{IIN14CWT}. In order to construct an example of non--zonal wavelets, it was necessary to modify the definition to the case of different analysis and reconstruction wavelet families~\cite{IIN15DW}. For the purposes of the present paper we need slightly different assumptions, similar to those given in~\cite{HH09}. The reason for the modification is that we want to investigate the continuous frame property of the wavelet transform which was not ensured in the case of different analysis and reconstruction families. Namely, it can be seen from the proof of Theorem~\ref{thm:energy_conservation} that a change of the exponent by~$\rho$ would break the frame property, whereas in the case of a bilinear wavelet according to~\cite[Definition~5.1]{IIN15DW} a change of the exponent by~$\rho$ in one of the families can be compensated by a proper change of the exponent by~$\rho$ in the second one.

\begin{df}\label{def:wavelet} Let $\alpha:\mathbb R_+\to\mathbb R_+$ be a weight function. The family $\{\Psi_\rho\}_{\rho\in\mathbb R_+}\subseteq\mathcal L^2(\mathbb S^n)$ is called a wavelet (family) of order~$m$ if it satisfies
\begin{equation}\label{eq:cond_energy_conservation}
A\leq\beta(l)\leq B
\end{equation}
for some positive constants~$A$ and~$B$ independent of $l\in\mathbb{N}_0$, $l>m$, and $\beta(l)=0$ for $l\in\mathbb N_0$, $l\leq m$, where $\beta(l)$ is defined by
$$
\beta(l):=\frac{\sum_{\kappa=1}^{N(n,l)}\int_0^\infty|a_l^\kappa(\Psi_\rho)|^2\,\alpha(\rho)d\rho}{N(n,l)}.
$$
\end{df}

For the convenience of the reader, we call $A$ and $B$ the \emph{wavelet family bounds} of the wavelet family~$\{\Psi_\rho\}$.

For such a wavelet one can perform a wavelet transform according to the following definition (cf. \cite[Definition~3.2]{IIN14CWT}).

\begin{df}\label{def:bilinear_wt} Let $\{\Psi_\rho\}_{\rho\in\mathbb R_+}$ be a wavelet family. Then, the spherical wavelet transform
$$
\mathcal W_\Psi\colon\mathcal L^2(\mathcal S^n)\to\mathcal L^2(\mathbb R_+\times SO(n+1))
$$
is defined by
$$
\mathcal W_\Psi f(\rho,\Upsilon)=\frac{1}{\Sigma_n}\int_{\mathcal S^n}\overline{\Psi_\rho(\Upsilon^{-1}x)}\,f(x)\,d\sigma_n(x).
$$
\end{df}

For the wavelet transform defined in this way a kind of energy conservation can be proven. In~\cite{HH09}, for the admissibility of the wavelet family the authors require the wavelet transform to be a frame, in our notation:
$$
A\|f\|^2\leq\int_0^\infty\int_{SO(3)}|\mathcal W_\Psi f(\rho,\Upsilon)|^2\,d\nu(\Upsilon)\,d\alpha(\rho)\leq B\|f\|^2.
$$
(It is worth noting that it is not exactly the same as boundedness of the wavelet transform or boundedness of the wavelet synthesis required in~\cite{mH96}, since the wavelet synthesis in~\cite{HH09} is not the inverse transform to the wavelet transform.) It is pointed out that this is equivalent to a condition on the Fourier coefficients of the wavelet family.  The present approach goes in another direction, i.e., we require some constraints on the wavelet family and prove that the frame property is satisfied.

\begin{thm}\label{thm:WT_frame_property} Let $\{\Psi_\rho\}$ be a wavelet family of order~$m$. Then for any $f\in\mathcal L^2(\mathcal S^n)$ with $m$ vanishing moments (i.e., such that $a_l^k(f)=0$ for $l=0,1,\dots,m$ and $k\in\mathcal M_{n-1}(l)$) we have
$$
A\|f\|^2\leq\|\mathcal W_\Psi f\|^2\leq B\|f\|^2,
$$
i.e, the set $\{\Psi_\rho(\Upsilon^{-1}\circ),\,\rho\in(0,\infty),\Upsilon\in SO(n+1)\}$ is a frame for $\mathcal L^2(\mathcal S^n)\setminus\bigcup_{l=0}^m\mathcal H_l(\mathcal S^n)$.
\end{thm}

\begin{bfseries}Proof. \end{bfseries}In the same manner as in the proof of \cite[Theorem~3.3]{IIN14CWT} we get
\begin{equation}\label{eq:WT_squared}\begin{split}
&\|\mathcal W_\Psi f\|^2\\
&=\int_{\mathcal S^n}\!\int_{\mathcal S^n}\!\int_0^\infty
   \frac{\sum_{l,k}\bigl|a_l^k(\Psi_\rho)\bigr|^2}{N(n,l)}\,\mathcal K_l^\lambda(x\cdot y)\,\alpha(\rho)\,d\rho\,\overline{f(x)}\,d\sigma_n(x)\,f(y)\,d\sigma_n(y).
\end{split}\end{equation}
By $\mathcal L^2$ integrability of~$\Psi_\rho$ and~$f$, and mutual orthogonality of~$f_l$'s, we thus obtain
\begin{equation}\label{eq:energy_conservation}
\|\mathcal W_\Psi f\|^2=\sum_{l=m+1}^\infty\beta(l)\,\int_{\mathcal S^n}\overline{f_l(y)}\,f_l(y)\,d\sigma_n(y).
\end{equation}
The claim follows from non--negativity of the integrals on the right--hand--side of~\eqref{eq:energy_conservation}, orthogonality of~$f_l$'s, inequalities~\eqref{eq:cond_energy_conservation}, and completeness of $\bigcup_{l=0}^\infty\mathcal H_l$ in $\mathcal L^2(\mathcal S^n)$.\hfill$\Box$\\

\label{rem:zonal_wavelets}For zonal wavelets appearing in the literature (cf., e.g., \cite{FGS-book,GM09a,IIN14PW}), condition~\eqref{eq:cond_energy_conservation} is satisfied with $A=B=1$, and the wavelet transform is an isometry. To the author's best knowledge, directional Poisson wavelets~\cite{HH09,IIN15DW} are the only so far constructed nonzonal wavelet family. In the theorem below we show that the method used for defining directional Poisson wavelets, i.e., directional derivation, can be applied to other zonal wavelet families and yields wavelets satisfying conditions of Definition~\ref{def:wavelet}. For the rest of the paper $\Upsilon_{\widehat\varsigma,\Theta}$ denotes the rotation of~$\mathbb R^{n+1}$ in the plane $(\widehat e,\widehat\varsigma)$, where $\widehat\varsigma$ is a rotation axis parallel to the tangent plane at~$\widehat e$, with rotation angle~$\Theta$, and set
$$
\Psi_\rho^{\widehat\varsigma,d}(x)=\rho^{ad/(\nu b)}\left.\frac{\partial^d}{\partial\Theta^d}\,\Psi_\rho(\Upsilon_{\widehat\varsigma,\Theta}x)\right|_{\Theta=0}.
$$

\begin{thm}\label{thm:energy_conservation} Let $\{\Psi_\rho\}\subset\mathcal L^2(\mathcal S^n)$ be a zonal family of functions satisfying
\begin{equation}\label{eq:cond_zonal_wavelets}
\widehat \Psi_\rho(l)=\left(\rho^a\,[q_\gamma(l)]^b\right)^c\,e^{-\rho^a\,[q_\gamma(l)]^b}\cdot\frac{l+\lambda}{\lambda}\qquad\text{for }l\in\mathbb N_0,
\end{equation}
where~$q_\gamma$ is a polynomial of degree~$\gamma$, strictly positive for a positive~$l$, and $a$, $b$, $c$ --- some positive constants. Then for any $d\in\mathbb N$, $\{\Psi_\rho^{\widehat\varsigma,d}\}$ is a wavelet of order~$0$ according to Definition~\ref{def:wavelet} with \mbox{$\alpha(\rho)=\frac{1}{\rho}$}.
\end{thm}

\begin{bfseries}Proof. \end{bfseries}Without loss of generality, we can assume~$\widehat\varsigma$ to be the $x_2$--axis. For the sake of simplicity set $\widetilde\rho=\rho^{a/(\gamma b)}$ and $\widetilde\Psi_{\widetilde\rho}=\Psi_\rho$. Note that $\frac{d\widetilde\rho}{\widetilde\rho}=\frac{d\rho}{\rho}$.
 In this case, according to \cite[Theorem~4.3]{IIN15DW} and \cite[Corollary~4.6]{IIN15DW}, the functions $\widetilde\Psi_{\widetilde\rho}^{[d]}=\widetilde\Psi_{\widetilde\rho}^{x_2,d}$ are given by
$$
\widetilde\Psi_{\widetilde\rho}^{[d]}
   =\widetilde\rho^d\,\sum_{l=0}^\infty\sum_{j=0}^{\left[\frac{d}{2}\right]}\widetilde a_l^{2j+d_{\text{mod}2}}(\widetilde\Psi_{\widetilde\rho}^{[d]})Y_l^{(2j+d_{\text{mod}2},0,\dots,0)}
$$
for $n\geq3$ or
$$
\widetilde\Psi_{\widetilde\rho}^{[d]}=\widetilde\rho^d\,\sum_{l=0}^\infty\sum_{j=0}^{\left[\frac{d}{2}\right]}\widetilde a_l^{2j+d_{\text{mod}2}}(\widetilde\Psi_{\widetilde\rho}^{[d]})
   \left(Y_l^{2j+d_{\text{mod}2}}+Y_l^{-(2j+d_{\text{mod}2})}\right)
$$
for $n=2$ with coefficients~$\widetilde a_l^j(\widetilde\Psi_{\widetilde\rho}^{[d]})$ of the form
$$
\widetilde a_l^j(\widetilde\Psi_{\widetilde\rho}^{[d]})
   =\left(\prod_{\iota=0}^{j-1}\beta_{l,\iota}\right)\cdot p_{d,j}\,(\beta_{l,0}^2,\beta_{l,1}^2,\dots,\beta_{l,\frac{d-j}{2}}^2)\cdot a_l^0(\tilde f_{\widetilde\rho})
$$
if $d$, $j$ are of the same parity, where~$p_{d,j}$ is a polynomial of degree~$\frac{d-j}{2}$, and
$$
\widetilde a_l^j(\widetilde\Psi_{\widetilde\rho}^{[d]})=0
$$
otherwise. Coefficients $\beta_{l,\iota}$ are given by
$$
\beta_{l,\iota}=\sqrt{\frac{(\iota +1)\,(2\lambda+\iota -1)}{(2\lambda+2\iota -1)\,(2\lambda+2\iota +1)}\,[l(2\lambda+l)-\iota (2\lambda+\iota )]},
$$
for $k_1=0,1,\dots,l$. Thus, $|\widetilde a_l^j(\widetilde\Psi_{\widetilde\rho}^{[d]})|^2$ is a polynomial in $l(2\lambda+l)$ of degree $j+2\cdot\frac{d-j}{2}=d$, multiplied by~$|a_l^0(\tilde f_{\widetilde\rho})|^2$ if $2|(d-j)$, and~$0$ otherwise. Consequently, by~\eqref{eq:cond_zonal_wavelets} and~\eqref{eq:hatfl_vs_al0f},
\begin{equation}\label{eq:beta_l_directional_wv}\begin{split}
\sum_{\kappa=1}^{N(n,l)}&\int_0^\infty|a_l^\kappa(\widetilde\Psi_{\widetilde\rho}^{[d]})|^2\,\frac{d\widetilde\rho}{\widetilde\rho}
   =\sum_{j=0}^d\int_0^\infty\widetilde\rho^{2d}\,|\widetilde a_l^j(\widetilde\Psi_{\widetilde\rho}^{[d]})|^2\,\frac{d\widetilde\rho}{\widetilde\rho}\\
&=P_{2d}(l)\cdot\left(\frac{\lambda+l}{\lambda\,A_l^0}\right)^2\int_0^\infty\widetilde\rho^{2d}\left[\widetilde\rho^\gamma q_\gamma(l)\right]^{2bc}
   e^{-2\left[\widetilde\rho^\gamma q_\gamma(l)\right]^b}\,\frac{d\widetilde\rho}{\widetilde\rho}
\end{split}\end{equation}
where~$P_{2d}$ is a polynomial of degree~$2d$. As a sum of squared modules, $P_{2d}$ is nonnegative. Suppose, $P_{2d}(l)$ is equal to~$0$ for some~$l>0$. This means
$$
(Y_l^0)^{(d)}:=\left.\frac{\partial^d}{\partial\Theta^d}\,Y_l^0(\Upsilon_{x_2,\Theta}x)\right|_{\Theta=0}=0.
$$
Thus,
$$
(Y_l^0)^{(d-1)}=\text{const.},
$$
i.e., either $(Y_l^0)^{(d-1)}$ is a multiple of~$Y_0^0$, which is a contradiction to $(Y_l^0)^{(d)}\in\mathcal H_l$, $l>0$, or it is equal to~$0$. By induction we obtain $Y_l^0=0$, a contradiction to $l>0$. The assumption $P_{2d}(l)=0$ for some $l>0$ was false, i.e., $P_{2d}$ is strictly positive for positive~$l$. On the other hand, it is easy to see that $P_{2d}(0)=0$ for $d\geq1$ since it corresponds to a derivative of a constant function.

Since $\left(\frac{\lambda+l}{\lambda\,A_l^0}\right)^2$ is equal to~$N(n,l)$, from~\eqref{eq:beta_l_directional_wv} we obtain
$$
\beta(l)=\frac{P_{2d}(l)\,\Gamma\left(2c+\frac{2d}{\gamma b}\right)}{4^{d/(\gamma b)+c}\,\gamma b\,[q_\gamma(l)]^{2d/\gamma}}.
$$
Both $q_\gamma(l)$ and $P_{2d}(l)$ are positive for positive~$l$, and their behavior in infinity is governed by~$l^{2d}$. Therefore, $\beta(l)$ is positive for $l>0$ and $\lim_{l\to\infty}\beta(l)$ exists. Thus, \eqref{eq:cond_energy_conservation} is satisfied for some positive constants $A$, $B$, independent of~$l\geq1$ and it vanishes for $l=0$.\hfill$\Box$

\begin{cor}\label{cor:directional_wavelets}Let a zonal family $\{\Psi_\rho\}$ of functions satisfying~\eqref{eq:cond_zonal_wavelets} be given, and let
$$
\Psi_\rho^{\widehat\varsigma,d}(x)=\rho^{ad/(\nu b)}\left.\frac{\partial^d}{\partial\Theta^d}\,\Psi_\rho(\Upsilon_{\widehat\varsigma,\Theta}x)\right|_{\Theta=0}.
$$
Then, $\{\Psi_\rho^{\widehat\varrho,d}(\Upsilon^{-1}\circ),\,\rho\in(0,\infty),\,\Upsilon\in SO(n+1)\}$ is a frame for $\mathcal L^2(\mathcal S^n)\setminus\mathcal H_0(\mathcal S^n)$.
\end{cor}

Corollary~\ref{cor:directional_wavelets} applies to all of the known zonal wavelet families, i.e., Abel--Poisson wavelets and Gauss--Weierstrass wavelets~\cite{FGS-book}, Poisson wavelets~\cite{IIN14PW}, as well as zonal diffusive wavelets characterized in \cite[Theorem~4.2.6]{sE11}. This means that the abandonment of the requirement of an exact reconstruction widens the set of wavelets significantly when compared to that investigated in~\cite{IIN15DW}. Another advantage of the present approach with respect to that from~\cite{IIN15DW} is that no artificial wavelet family is introduced, and one can use simple derivatives of directional wavelet rotations instead of their linear combinations.

\section{Semi--continuous wavelet frames}\label{sec:semi_cont_frames}
The goal of this paper is to prove that wavelets constructed as in Theorem~\ref{thm:energy_conservation} possess discrete frames. First, we discretize the scale parameter~$\rho$ -- in a similar manner as in~\cite[Section~3]{IIN15WF}.

\begin{thm}
Let~$\{\Psi_\rho\}\subset\mathcal L^2(\mathcal S^n)$ be a wavelet family as in Theorem~\ref{thm:energy_conservation} with wavelet family bounds~$A$ and~$B$. Then, for any $\epsilon>0$ there exist constants~$\mathfrak a_0$ and $X$ such that for any sequence~$\mathcal{R}=(\rho_j)_{j\in\mathbb{N}_0}$ with $\rho_0\geq\mathfrak a_0$ and $1<\rho_j/\rho_{j+1}<X$ the family \mbox{$\{\Psi_{\rho_j}(\Upsilon^{-1}\circ),\,\rho_j\in\mathcal{R},\,\Upsilon\in SO(n+1)\}$} is a semi-continuous frame for~$\mathcal{L}^2(\mathcal S^n)$, satisfying the frame condition~\eqref{eq:framedef_AB} with the frame bounds~$A(1-\epsilon)$ and~$B(1+\epsilon)$.
\end{thm}

\begin{bfseries}Proof.\end{bfseries} Replace~$\rho$ by~$\widetilde\rho=\rho^{a/(\nu b)}$. The discretization of the scale parameter leads to the discretization of the integral of $\sum_\kappa|a_l^\kappa|^2$ in~\eqref{eq:WT_squared}, and in the present case it is equivalent to the discretization of the integral in~\eqref{eq:beta_l_directional_wv}. A proper choice of weights and parameters~$\mathfrak a_0$ and~$X$ ensures that the sum of the series is arbitrarily close (i.e., $<\epsilon$) to the value of the integral, compare the proof of \cite[Theorem~5]{IH10} (with $\widetilde\rho^{2d}\left[\widetilde\rho^\nu q_\nu(l)\right]^{bc}e^{-\left[\widetilde\rho^\nu q_\nu(l)\right]^b}$ substituted for~$\gamma_d(\widetilde\rho)$).\hfill$\Box$

\section{Discretization of rotation parameter}\label{sec:discrete_frames}

In order to discretize the rotation parameter $\Upsilon\in SO(n+1)$, we represent it as
$$
\Upsilon=\Upsilon^{1}\,\Upsilon^{2}\dots\Upsilon^{n}
$$
with
$$
\Upsilon^J(x^J)=\Upsilon_1(\theta_1^{J})\,\Upsilon_2(\theta_2^{J})\dots\Upsilon_{J-1}(\theta_{J-1}^{J})\,\Upsilon_J(\varphi^{J}),\qquad J=1,\,2,\dots,\,n,
$$
for
\begin{align*}
&x^1=\varphi^1\in\mathcal S^1,\quad x^2=(\theta_1^2,\varphi^2)\in\mathcal S^2,\quad\dots,\\
&x^{n-1}=(\theta_1^{n-1},\dots,\theta_{n-2}^{n-1},\varphi^{n-1})\in\mathcal S^{n-1},\quad x^n=(\theta_1^{n},\dots,\theta_{n-1}^{n},\varphi^{n})\in\mathcal S^n.
\end{align*}
where~$\Upsilon_\iota(\theta)$ is the rotation in the plane $(x_\iota,x_{\iota+1})$ with the rotation angle~$\theta$, and $\theta_\iota^{J}$ and $\varphi^{J}$, $J=1,2\dots,n$, $\iota=1,2,\dots,J-1$, are the Euler angles of~$\Upsilon$, according to \cite[Theorem~IX.1.1]{Vilenkin}. We interpret $SO(n+1)$ as the Cartesian product of~$\mathcal S^J$, $J=1,\,2,\dots,\,n$, and discretization will be performed separately on each~$\mathcal S^J$.

\begin{df}\label{df:grid}Let $n$ be given. We say $\Lambda$ is a grid of type $(\delta_n,\delta_{n-1},\dots,\delta_1)$ if it is a discrete measurable set of rotations in~$SO(n+1)$, constructed iteratively in the following way. There is a measurable partition $\mathcal P_n=\{\mathcal O_{\alpha_n}^n:\,\alpha_n=1,\dots,K_n\}$ of~$\mathcal S^n$ into simply connected sets such that the diameter of each set (measured in geodesic distance) is not larger than~$\delta_n$. Each of these sets contains exactly one point~$x_{\alpha_n}^n$, $\alpha_n=1,\dots,K_n$. For $J<n$ let $(\alpha_n,\alpha_{n-1},\dots,\alpha_{J+1})$ be a fixed multi-index and
$$
\mathcal P_J=\mathcal P_j(\alpha_n,\alpha_{n-1},\dots,\alpha_{J+1})=\{\mathcal O_{\alpha_J}^J:\,\alpha_J=1,\dots,K_J\}
$$
a  measurable partition of~$\mathcal S^J$ into $K_J=K_J(\alpha_n,\alpha_{n-1},\dots,\alpha_{J+1})$ simply connected sets of diameter not larger than~$\delta_J$, each of them containing exactly one point $x_{\alpha_J}^J=x_{(\alpha_n,\dots,\alpha_J)}^J$. Then, $\Lambda$ is the set of rotations given by
$$
\Upsilon_{(\alpha_n,\dots,\alpha_1)}=\Upsilon^{1}(x_{\alpha_1}^1)\,\Upsilon^{2}(x_{\alpha_2}^2)\dots\Upsilon^{n}(x_{\alpha_n}^n)
$$
with the measure
$$
\lambda(\Upsilon_{(\alpha_n,\dots,\alpha_1)}):=\prod_{J=1}^n\lambda^J\left(x_{\alpha_J}^J\right),
\qquad\lambda^J\left(x_{\alpha_J}^J\right):=\sigma_J\left(\mathcal O_{\alpha_J}^J\right).
$$
\end{df}

\begin{thm} Let $\Psi_\rho$ be a $\mathcal C^1$--wavelet family with the property that $\{\Psi_{\rho_j,x},\,j\in\mathbb N_0,\,x\in\mathcal S^n\}$ is a semi--continuous frame. Then, for each $j\in\mathbf N_0$ there exist sequences $(\delta_n^j,\delta_{n-1}^j,\dots,\delta_1^j)$ such that
$$
\{\Psi_{\rho_j}(\Upsilon_{(\alpha_n^j,\dots,\alpha_1^j)}^{-1}\circ),\,j\in\mathbf N_0,\,\Upsilon_{(\alpha_n^j,\dots,\alpha_1^j)}\in\Lambda^j\}
$$
is a frame for~$\mathcal L^2(\mathcal S^n)$, provided that $\Lambda^j$ is a grid of type $(\delta_n^j,\delta_{n-1}^j,\dots,\delta_1^j)$.
\end{thm}

\begin{bfseries}Proof. \end{bfseries}Notation: $\mathfrak c$, resp. $\mathfrak c_J$, is a positive constant that may change its value from line to line.\\
Let~$\rho_j$ be fixed. Then,
\begin{align*}
\int_{SO(n+1)}&|\mathcal W_\Psi(\rho_j,\Upsilon)|^2\,d\nu(\Upsilon)
   =\frac{1}{\Sigma_n^2}\int_{SO(n+1)}\left|\int_{\mathcal S^n}\overline{\Psi_{\rho_j}(\Upsilon^{-1}x)}\,f(x)\,d\sigma_n(x)\right|^2d\nu(\Upsilon)\\
=&\frac{1}{\Sigma_n^2\prod_{J=1}^n\Sigma_J}\int_{\mathcal S^n}\int_{\mathcal S^{n-1}}\dots\int_{\mathcal S^2}\int_{\mathcal S^1}\\
&\left[\int_{\mathcal S^n}\overline{\Psi_{\rho_j}\left((\Upsilon^{n})^{-1}\,(\Upsilon^{n-1})^{-1}\dots(\Upsilon^{2})^{-1}\,(\Upsilon^{1})^{-1}\,x\right)}
   f(x)\,d\sigma_n(x)\right]^2\\
&d\sigma_1(x^1)\,d\sigma_2(x^2)\dots d\sigma_{n-1}(x^{n-1})\,d\sigma_n(x^n).
\end{align*}
We will prove by induction that the integral can be discretized and the discretization error is less than $2^{-j-1}\,\delta\,\|f\|_2^2$ for a prescribed~$\delta$, provided that the constants $\delta_n^j,\dots,\delta_1^j$ are small enough. Set
$$
\mathcal Y^J=\Upsilon^J\,\Upsilon^{J+1}\dots\Upsilon^{n-1}\,\Upsilon^n,\qquad\mathcal Y^{n+1}=Id
$$
and
$$
\mathcal X^J=\Upsilon^1\,\Upsilon^2\dots\Upsilon^{J-1}\,\Upsilon^J,\qquad\mathcal X^0=Id,
$$
then for each $J=1,\,2,\dots,n$
$$
\Upsilon=\mathcal X^{J-1}\Upsilon^J\mathcal Y^{J+1}.
$$
For $J=n$ consider the integral
\begin{align*}
I_J&:=\int_{\mathcal S^J}\int_{\mathcal S^{J-1}}\dots\int_{\mathcal S^1}\left[\int_{\mathcal S^n}
   \overline{\Psi_{\rho_j}(\Upsilon^{-1}x)}\,f(x)\,d\sigma_n(x)\right]^2d\sigma_1(x^1)\dots d\sigma_J\!\left(x^J\right)\\
&=\int_{\mathcal S^J}\int_{\mathcal S^{J-1}}\dots\int_{\mathcal S^1}
   \left[\int_{\mathcal S^n}\overline{\Psi_{\rho_j}\left((\Upsilon^J)^{-1}\,\left(\mathcal X^{J-1}\right)^{-1}\,x\right)}\,f(\mathcal Y^{J+1}\,x)\,d\sigma_n(x)\right]^2\\
&\qquad d\sigma_1(x^1)\dots d\sigma_J\left(x^J\right)
\end{align*}
and set
\begin{equation*}\begin{split}
&E^J=E_J(\alpha_n,\dots,\alpha_{J+1}):=I_J\\
&-\sum_{\alpha_J=1}^{K_J(\mathcal Y^{J+1})}\int_{\mathcal S^{J-1}}\dots\int_{\mathcal S^1}
\left[\int_{\mathcal S^n}
   \overline{\Psi_{\rho_j}\left(\left(\Upsilon^J(x_{\alpha_J}^J)\right)^{-1}\left(\mathcal X^{J-1}\right)^{-1}x\right)}\,f(\mathcal Y^{J+1}\,x)\,
   d\sigma_n(x)\right]^2\\
&\qquad d\sigma_1(x^1)\dots d\sigma_{J-1}\left(x^{J-1}\right)\lambda^J(x_{\alpha_J}^J).
\end{split}\end{equation*}
Then,
\begin{equation}\label{eq:est_E1}\begin{split}
|E^J|\leq&\sum_{\alpha_J=1}^{K_J}\int_{\mathcal O_{\alpha_J}^J}
   \int_{\mathcal S^{J-1}}\dots\int_{\mathcal S^1}\\
&\left|\left[\int_{\mathcal S^n}\overline{\Psi_{\rho_j}\left(\left(\Upsilon^J(x_{\alpha_J}^J)\right)^{-1}\left(\mathcal X^{J-1}\right)^{-1}x\right)}f(\mathcal Y^{J+1}\,x)\,d\sigma_n(x)\right]^2\right.\\
&-\left.\left[\int_{\mathcal S^n}\overline{\Psi_{\rho_j}\left(\left(\Upsilon^J(x^J)\right)^{-1}\left(\mathcal X^{J-1}\right)^{-1}x\right)}
   f(\mathcal Y^{J+1}\,x)\,d\sigma_n(x)\right]^2\right|\\
&d\sigma_1(x^1)\dots d\sigma_{J-1}\left(x^{J-1}\right)d\sigma_J(x^J).
\end{split}\end{equation}
For each~$\alpha_J$, the integrand can be estimated by
\begin{equation}\label{eq:est_e_alpha_J}\begin{split}
&e_{\alpha_J}^J=e_{\alpha_J}^J(\alpha_n,\dots,\alpha_{J+1},\mathcal X^{J-1})\\
&:=\sup_{x^J\in\mathcal O_{\alpha_J}^J}\left|\nabla^\ast
   \left[\int_{\mathcal S^n}\overline{\Psi_{\rho_j}\left(\left(\Upsilon^J(x^J)\right)^{-1}\left(\mathcal X^{J-1}\right)^{-1}x\right)}
   f(\mathcal Y^{J+1}\,x)\,d\sigma_n(x)\right]^2\right|\cdot\delta\left(\mathcal O_{\alpha_J}^J\right)\\
&\leq2\int_{\mathcal S^n}\sup_{y\in\mathcal S^n}\left|\Psi_{\rho_j}(y)\right|\left|f(\mathcal Y^{J+1}\,x)\right|\,d\sigma_n(x)\\
&\cdot\int_{\mathcal S^n}\sup_{y\in\mathcal S^n}\left|\nabla^\ast\Psi_{\rho_j}(y)\right|
   \left|f(\mathcal Y^{J+1}\,x)\right|\,d\sigma_n(x)\cdot\delta_J^j,
\end{split}\end{equation}
and, thus, by H\"older inequality
\begin{equation}\label{eq:est_ealpha1}
e_{\alpha_J}^J\leq\mathfrak c\cdot\|\sup\Psi_{\rho_j}\|_2\cdot\|f\|_2\cdot\|\sup\nabla^\ast\Psi_{\rho_j}\|_2\cdot\|f\|_2
   \cdot\delta_J^j.
\end{equation}
Consequently, since $\sum_{\alpha_J=1}^{K_J}\sigma_J\left(\mathcal O_{\alpha_J}^J\right)=\Sigma_J$, from~\eqref{eq:est_E1} and~\eqref{eq:est_ealpha1} we have
\begin{equation}\label{eq:estEJ}
|E^J|\leq\mathfrak c\cdot\|\Psi_{\rho_j}\|_\infty\cdot\|\nabla^\ast\Psi_{\rho_j}\|_\infty\cdot\delta_J^j\cdot\|f\|_2^2=:E_J.
\end{equation}
Now, suppose, for some $J<n$ the following holds
\begin{align*}
&\int_{SO(n+1)}|\mathcal W_\Psi(\rho_j,\Upsilon)|^2\,d\nu(\Upsilon)\\
&=\frac{1}{\Sigma_n^2\prod_{J=1}^n\Sigma_J}\sum_{\alpha_n=1}^{K_n}\dots
   \Biggl(\sum_{\alpha_{J+2}=1}^{K_{J+2}(\mathcal Y^{J+3})}\Biggl(\sum_{\alpha_{J+1}=1}^{K_{J+1}(\mathcal Y^{J+2})}
   \int_{\mathcal S^J}\dots\int_{\mathcal S^1}\\
&\left[\int_{\mathcal S^n}\Psi_{\rho_j}\left(\left(\Upsilon^J(x^J)\right)^{-1}\dots\left(\Upsilon^1(x^1)\right)^{-1}x\right)
\cdot f\left(\Upsilon^{J+1}(x_{\alpha_{J+1}}^{J+1})\dots\Upsilon^n(x_{\alpha_n}^n)\,x\right)\,d\sigma_n(x)\right]^2\\
&d\sigma_1(x^1)\dots d\sigma_J(x^J)\,\lambda^{J+1}(x_{\alpha_{J+1}}^{J+1})+E^{J+1}\Biggr)\lambda^{J+2}(x_{\alpha_{J+2}}^{J+2})+E^{J+2}\Biggr)
   \dots\lambda^n(x_{\alpha_n}^n)+E^n
\end{align*}
with~$|E^{J+1}|$ uniformly bounded by~$E_{J+1}$. For each multi--index $(\alpha_n,\dots,\alpha_{J+1})$ the integral
$$
I_J=\int_{\mathcal S^J}\int_{\mathcal S^{J-1}}\dots\int_{\mathcal S^1}\Bigl[\dots\Bigr]^2\,d\sigma_1(x^1)\dots d\sigma_{J-1}(x^{J-1})\,d\sigma_J(x^J)
$$
can be replaced by
$$
\sum_{\alpha_J=1}^{K_J}\int_{\mathcal S^{J-1}}\dots\int_{\mathcal S^1}\Bigl[\dots\Bigr]^2\,
   d\sigma_1(x^1)\dots d\sigma_{J-1}(x^{J-1})\,\lambda^J(x_{\alpha_J}^J).
$$
The discretization error is equal to~$E^J$ and in its absolute value bounded by~$E_J$. The statement can be proven along the same lines as for $J=n$. Thus, by induction,
\begin{align*}
\int_{SO(n+1)}&|\mathcal W_\Psi(\rho_j,\Upsilon)|^2\,d\nu(\Upsilon)=\frac{1}{\Sigma_n^2\prod_{J=1}^n\Sigma_J}\\
\cdot&\sum_{\alpha_n=1}^{K_n}\dots\sum_{\alpha_1=1}^{K_1}
   \left[\int_{\mathcal S^n}\Psi_{\rho_j}(x)\cdot f\left(\Upsilon^1(x_{\alpha_1}^1)\dots\Upsilon^n(x_{\alpha_n}^n)\,x\right)\,d\sigma_n(x)\right]^2\\
&\lambda^1(x_{\alpha_1}^1)\dots\lambda^n(x_{\alpha_n}^n)+\mathcal E^n
\end{align*}
with
$$
\mathcal E^n=\sum_{\alpha_n=1}^{K_n}\left(\dots\left(\sum_{\alpha_3=1}^{K_3}
   \left(\sum_{\alpha_2=1}^{K_2}E^1\,\lambda^2(x_{\alpha_2}^2)+E^2\right)\lambda^3(x_{\alpha_3}^3)+E^3\right)\dots\right)\dots+E^n,
$$
Since $|E^J|$ is bounded by~$E_J$, and $\sum_{\alpha_J=1}^{K_J}\lambda^J(x_{\alpha_J}^J)=\Sigma_J$,
$$
|\mathcal E^n|\leq\left(\dots\left((E_1\Sigma_2+E_2)\,\Sigma_3+E_3\right)\dots\right)+E_n.
$$
Consequently, by~\eqref{eq:estEJ},
$$
|\mathcal E^n|\leq\sum_{J=1}^n\mathfrak c_J\cdot E_J=\sum_{J=1}^n\mathfrak c_J\cdot\delta_J^j\cdot\|\Psi_{\rho_j}\|_\infty\cdot\|\nabla^\ast\Psi_{\rho_j}\|_\infty\cdot\|f\|_2^2.
$$
A proper choice of $\delta_J^j$'s ensures that
$$
|\mathcal E^n|\leq2^{-j-1}\,\delta\,\|f\|_2^2.
$$
The assertion follows from Theorem~\ref{thm:frames_discretization}.\hfill$\Box$\\

Although one could think about a finer error estimation in~\eqref{eq:est_e_alpha_J}, the author was not able to find a satisfactory method applicable in this general setting. It is very probable that for a concrete wavelet family, the error growth can be controlled by some power of the scale parameter~$\rho_j$ (similarly as in~\cite{IIN15WF}), and an upper bound of the grid density can be given. Obviously, this is possible also using inequality~\eqref{eq:est_e_alpha_J}, but then the obtained sufficient grid density would be far from being optimal, and the author prefers to avoid complicated calculations.

\end{document}